\documentclass[graybox]{svmult}
\usepackage{graphicx}       % standard LaTeX graphics tool
\usepackage[utf8]{inputenc}
\usepackage[T1]{fontenc}
\usepackage{microtype} % good font tricks

\usepackage{cite}
\usepackage{array,colortbl}

\usepackage{amsmath,amsfonts,amssymb,bm,mathtools,mathrsfs}
\DeclareFontFamily{U}{mathx}{\hyphenchar\font45}
\DeclareFontShape{U}{mathx}{m}{n}{
      <5> <6> <7> <8> <9> <10>
      <10.95> <12> <14.4> <17.28> <20.74> <24.88>
      mathx10
      }{}
\DeclareSymbolFont{mathx}{U}{mathx}{m}{n}
\DeclareFontSubstitution{U}{mathx}{m}{n}
\DeclareMathAccent{\widecheck}{0}{mathx}{"71}

\def\citep#1#2{\cite[{#1}]{#2}}

\newcommand{\RefEq}[1]{~\textup{(\ref{#1})}}
\newcommand{\RefEqTwo}[2]{~\textup{(\ref{#1})} and~\textup{(\ref{#2})}}

\newcommand{\RefSec}[1]{Section~\textup{\ref{#1}}}

\newcommand{\RefThm}[1]{Theorem~\textup{\ref{#1}}}

\newcommand{\RefFig}[1]{Figure~\textup{\ref{#1}}}

\usepackage{algorithm}
\usepackage{algorithmic}
\newcommand{\N}{{\mathbb{N}}} % natural numbers {1, 2, ...}
\newcommand{\R}{{\mathbb{R}}} % reals
\DeclareSymbolFont{bbold}{U}{bbold}{m}{n}
\DeclareSymbolFontAlphabet{\mathbbold}{bbold}

\newcommand{\calA}{{\mathcal{A}}}

\newcommand{\calC}{{\mathcal{C}}}
\newcommand{\calD}{{\mathcal{D}}}

\newcommand{\calP}{{\mathcal{P}}}

\newcommand{\calR}{{\mathcal{R}}}

\usepackage{type1cm}        % activate if the above 3 fonts are
\usepackage{makeidx}         % allows index generation
\usepackage{graphicx}        % standard LaTeX graphics tool
\usepackage{multicol}        % used for the two-column index
\usepackage[bottom]{footmisc}% places footnotes at page bottom

\usepackage{newtxtext}       % 
\usepackage[varvw]{newtxmath}       % selects Times Roman as basic font

\usepackage{dsfont}
\makeindex             % used for the subject index
\newcommand{\mathset}[2]{\left\{#1\middle\vert #2 \right\}}
\newcommand{\Exp}[1]{\mathbb{E}\left(#1\right)}

\newcommand{\mathseq}[2]{\left(#1\right)_{#2}}
\newcommand{\Prob}[1]{\mathbb{P}\left(#1\right)}
\newcommand{\card}[1]{{\left\vert #1 \right\vert}}
\newcommand{\leb}{\mathbb{\lambda}}

\newcommand*\diff{\mathop{}\!\mathrm{d}}

\DeclareMathOperator{\variance}{Var}

\begin{document}

\title*{New Bounds for the Extreme and the Star Discrepancy of Double-Infinite Matrices}
\titlerunning{New Bounds for the Discrepancy of Double-Infinite Matrices}

%\title*{New Bounds for the extreme Discrepancy and Star-Discrepancy of Double-Infinite Matrices}

\date{\today}

\author{Jasmin Fiedler and Michael Gnewuch and Christian Wei\ss{}}
\institute{Jasmin Fiedler \at Institut Naturwissenschaften, Ruhr West University of Applied Sciences, Duisburger~Str.~100, 45479 M\"ulheim an der Ruhr, \email{manuel.fiedler@hs-ruhrwest.de}
\and Michael Gnewuch \at Institut f\"ur Mathematik, Osnabr\"uck University, Albrechtstr.~28a, 49076 Osnabr\"uck, \email{michael.gnewuch@uni-osnabrueck.de} \and Christian Wei\ss{} \at Institut Naturwissenschaften, Ruhr West University of Applied Sciences, Duisburger~Str.~100, 45479 M\"ulheim an der Ruhr, \email{christian.weiss@hs-ruhrwest.de}}

\maketitle

%\author{Manuel Fiedler, Michael Gnewuch and Christian Weiß}
%\address{{\bf{Ruhr West University of Applied Sciences,}}\\ {{Department of Natural Sciences, Duisburger Str. 100,}}\\{{45479 M\"ulheim an der Ruhr, Germany}}}
%\email{manuel.fiedler@hs-ruhrwest.de}
%\maketitle

\abstract{ According to Aistleitner and Weimar, there exist two-dimensional (double) infinite matrices whose star-discrepancy $D_N^{*s}$ of the first $N$ rows and $s$ columns, interpreted as $N$ points in $[0,1]^s$, satisfies an inequality of the form 
$$D_N^{*s} \leq \sqrt{\alpha} \sqrt{A+B\frac{\ln(\log_2(N))}{s}}\sqrt{\frac{s}{N}}$$
with $\alpha = \zeta^{-1}(2) \approx 1.73, A=1165$ and $B=178$. These matrices are obtained by using i.i.d sequences, and the parameters $s$ and $N$ refer to the dimension and the sample size respectively. In this paper, we improve their result in two directions: First, we change the character of the equation so that the constant $A$ gets replaced by a value $A_s$ dependent on the dimension $s$ such that for $s>1$ we have $A_s<A$. Second, we generalize the result to the case of the (extreme) discrepancy. The paper is complemented by a section where we show numerical results for the dependence of the parameter $A_s$ on $s$.}%we generalize the result towards a wider class of sets, most notably the one describing the extreme discrepancy.}

\keywords{(extreme) discrepancy, star-discrepancy, random matrices, bracketing numbers}



\section{Introduction}
The task of high-dimensional integration occurs in different practical applications, with one of them being computational finance (see e.g. \cite{Gla03}, \cite{DHB20}, \cite{Pas94}, \cite{WN19}). Since it is often impossible to give analytical solutions of the integrals, it is necessary to use numerical integration instead. Then the\textit{ goodness} of the approximation is in general mainly governed by the star-discrepancy. 
\begin{definition} For two points $x,y\in [0,1]^s$, we write $x\leq y$ if the inequality holds component-wise and we set $[0,x):=\mathset{y\in [0,1)^s}{0\leq y< x}$ and $[x,y):=\mathset{z\in [0,1)^s}{x\leq z< y}$. Let $(x_1,\ldots,x_n)$ be a sequence of points in $[0,1)^s$. Then we define their (extreme) discrepancy as
\begin{equation}
    D_N^{s}(x_1,\ldots,x_N)=\sup_{a,b\in[0,1]^s}\left\vert\frac{1}{N}\sum_{k=1}^N\mathds{1}_{[a,b)}\left(x_k\right)-\leb([a,b))\right\vert,
\end{equation}
where $\mathds{1}_{[a,b)}$ is the characteristic function of the interval $[a,b)$ and $\leb$ is the $s$-dimensional Lebesgue measure. Similarly, we define the star-discrepancy as
\begin{equation}
    D_N^{*s}(x_1,\ldots,x_N)=\sup_{y\in[0,1]^s}\left\vert\frac{1}{N}\sum_{k=1}^N\mathds{1}_{[0,y)}\left(x_k\right)-\leb([0,y))\right\vert.
\end{equation}
\end{definition}
The extreme discrepancy is also known under different names such as unanchored discrepancy, see e.g. \cite[Introduction]{chen:panorama_of_discr_theory}. 
\\
The precise dependence of the approximation error on the star-discrepancy is given by the iconic Koksma-Hlawka inequality, see e.g. \cite[Theorem 2.11 and 2.12]{nieder:rand_num_quasi_monte_carlo}.
\begin{theorem}[Koksma-Hlawka inequality] Let $f$ be a function on $[0,1]^s$ with bounded variation $V(f)$ in the sense of Hardy and Krause . Then for any sequence $(x_1,\ldots,x_n)$ of points in $[0,1)^s$ the inequality
\begin{equation}
    \left\vert\frac{1}{n}\sum_{k=1}^n f(x_k)-\int_{[0,1]^s}f(x)\diff x\right\vert\leq V(f)\cdot D_n^{*s}(x_1,\ldots,x_n)
\end{equation}
holds and the inequality is sharp.
\end{theorem}
Classically, low-discrepancy point sets are used for numerical integration. These are sets which satisfy the conjectured (known as the \textit{great open problem} of discrepancy theory) optimal bound 
\[D_N^{*s}(x_1,\ldots,x_N) \leq c_s\frac{\log(N)^{s-1}}{N},\]
where $c_s$ is a constant only depending on the dimension $s$. In fact, the conjecture is trivially true for dimension $1$ and has been proven to be correct in dimension $2$ by \cite{Sch72}. However, these low-discrepancy sets (or more precisely the upper bound of their star-discrepancy) suffer from the so-called \textit{curse of dimensionality}, i.e. the star-discrepancy exponentially depends on the dimension. If the sample size is small in comparison to the dimension, they are therefore hardly of practical use. In the paper \cite{HNWW01}, an alternative approach was therefore suggested. Therein, it was theoretically shown that the smallest achievable star-discrepancy with explicitly given dependence on the number of points as well as on the dimension satisfies
\[D_N^{*s}(x_1,\ldots,x_N) \leq C \sqrt{\frac{s}{N}}\]
without giving an explicit value for $C$. The first explicit value for $C \approx 9.65$ was derived in \cite{Ais11}, while currently, the best known possible value is $C \approx 2.4968$ due to \cite{gnew:gen_faul_ineq}. In \cite{Dic07, DGKP08} corresponding bounds were derived for double infinite matrices. The best currently known bound is due to \cite{aistl:prob_star_discr_double_inf_mat}.
%These will be used in order to adapt a result from \cite{aistl:prob_star_discr_double_inf_mat}. In  the following was shown:
\begin{remark} In this paper, $\zeta:(1,\infty)\to\R$ denotes the Riemann Zeta function on the reals, which is strictly decreasing, so the inverse $\zeta^{-1}$ exists.
\end{remark}
\begin{theorem}[Aistleitner, Weimar, \cite{aistl:prob_star_discr_double_inf_mat}, Theorem 1] \label{fi_thm:aw}
Let $\alpha>\zeta^{-1}(2)\approx 1.73$ be arbitrarily fixed. Then with probability strictly larger than $1- (\zeta(\alpha)-1)^2 \geq 0$ the double-infinite matrix $\mathseq{X_{n,s}}{n,s\in\N}$ with all $X_{n,s}\in[0,1)$ independently uniformly distributed satisfies for all $S \in \N$ and every $N \geq 2$ we have
\begin{equation}\label{fi_eq:aistl_ineq_discrep}
    D_N^{*S}(\calP_{N,S})\leq\sqrt{\alpha} \sqrt{A+B\frac{\ln(\log_2(N))}{S}}\sqrt{\frac{S}{N}},
\end{equation}
where $A=1165$, $B=178$ and
\[
\calP_{N,S}=\mathseq{X_{n,s}}{n,s=1}^{N,S}=\begin{pmatrix}
X_{1,1}&X_{1,2}&\cdots &X_{1,S}\\
X_{2,1}&X_{2,2}&\cdots &X_{2,S}\\
\vdots&\vdots&\ddots&\vdots\\
X_{N,1}&X_{N,2}&\cdots &X_{N,S}
\end{pmatrix}.
\]
\end{theorem}
The proof in \cite{aistl:prob_star_discr_double_inf_mat} relied on using the so-called maximal Bernstein inequality for independent random variables which will also be a main tool for our paper.
\begin{theorem}[Einmahl, Mason, \cite{einmahl:univ_res_beh_inc_part_sums}, Lemma 2.2]\label{fi_thm:max_bernstein_ineq} Let $X_1,\ldots,X_n$ be independent random variables with  $\Exp{X_k}=0$, $\left\vert X_k\right\vert\leq M<\infty$ a.s. and $\variance(X_k)=\sigma_k^2<\infty$ for all $1\leq k\leq n$. Then for every $t>0$
\begin{equation}\label{fi_eq:max_bernstein_ineq}
\Prob{\max_{1\leq k\leq n}\left\vert\sum_{i=1}^k X_i\right\vert >t}\leq 2\exp\left(-\frac{t^2}{2\sum_{k=1}^n \sigma_k^2+\frac23 Mt}\right).
\end{equation}
\end{theorem}
By combining the maximal Bernstein inequality with recent results proven in \cite{gnew:gen_faul_ineq}, we are able to improve the numerical constants in the formulation of \RefThm{fi_thm:aw}.
\begin{theorem}\label{fi_thm:discrep_bound_star}Let $\mathseq{X_{n,s}}{n,s\in\N}$ be a sequence of independent uniformly distributed random variables in $[0,1)$. Furthermore, take $\alpha> 1$ and $\beta>1$ with \\$(\zeta(\alpha)-1)(\zeta(\beta)-1)<1$. Then we have
\begin{align*}
\Prob{D_N^{*s}(\calP_{N,s})\leq \sqrt{\alpha A_{s}+\beta B\frac{\ln (\log_2(N))}{s}}\sqrt{\frac{s}{N}}\mathrm{~for~all~}s,N\in\N}\\\
\hspace{5cm}>1-(\zeta(\alpha)-1)(\zeta(\beta)-1)
\end{align*}
where $B=178$ and $A_{2}\leq 942$, $A_s$ is dependent on $s$  and decreases as $s$ increases.
\end{theorem}
Note that for dimension $s=1$, the result is not of a big practical relevance, because in this case low-discrepancy sequences (see e.g. \cite{nieder:rand_num_quasi_monte_carlo} for a definition) yield better results than random  sequences 
%do not suffer the curse of dimensionality then 
and  would therefore be preferred anyhow. In contrast to \RefThm{fi_thm:aw}, the constant $A_s$ is chosen dependent on the dimension $s$ such that for $s>1$ we have $A_s<A$. In all practically relevant cases, our result is therefore an improvement in comparison to \cite{aistl:inverse_discr_inf_dim_inf_seq}. Detailed numerical results on the constant $A_s$ will be discussed in \RefSec{fi_cp:num_res}.\\[12pt]
In addition and even more importantly, the methods applied here also cover the case of the (extreme) discrepancy. It is known that
\[
D_N^s(x_1,\ldots,x_N)\leq 2^s D_N^{*s}(x_1,\ldots,x_N),
\]
see for example \cite[1.4 Observation]{matousek:geometric_discrepancy_illustrated}, but this approximation is not useful in determining bounds for the discrepancy that scale well with the dimension. Therefore, we are going to use a different approach in the proof of \RefThm{fi_thm:discrep_bound_not_star}.
\begin{theorem}\label{fi_thm:discrep_bound_not_star}Let $\mathseq{X_{n,s}}{n,s\in\N}$ be a sequence of independent, uniformly distributed random variables in $[0,1)$. Furthermore, take $\alpha> 1$ and $\beta>1$ with \\$(\zeta(\alpha)-1)(\zeta(\beta)-1)<1$. Then we have
\begin{align*}
\Prob{D_N^{s}(\calP_{N,s})\leq \sqrt{\alpha A_{s}+\beta B\frac{\ln (\log_2(N))}{s}}\sqrt{\frac{s}{N}}\mathrm{~for~all~}s,N\in\N}\\\
\hspace{5cm}>1-(\zeta(\alpha)-1)(\zeta(\beta)-1)
\end{align*}
where $B=178$ and $A_{s}\leq 2548$ is dependent on $s$  and decreases as $s$ increases.
\end{theorem}

The paper is organized as follows: In \RefSec{fi_cp:star_disc} we will apply \RefThm{fi_thm:max_bernstein_ineq} in order to derive bounds for the discrepancy and the star-discrepancy of double infinite matrices. Afterwards, in \RefSec{fi_cp:num_res}, we will calculate the numerical values for the implied constants of the discrepancy bounds.

\paragraph{Acknowledgements} The authors would like to thank the organizers of the MCQMC Workshop 2022 in Linz, where a talk on a preliminary version of this article was given by the first-named author. Moreover, we would like to thank Bence Borda for asking a question which led to the generalization of our result from the star-discrepancy to the (extreme) discrepancy. Furthermore, we would like to thank the anonymous referee for their helpful comments. 

%\subsection{Latin hypercube samples} \label{subsec:exa}
%One prominent example for $\gamma$-negatively dependent samples is that of Latin hypercube samples.\textbf{CW: Quelle f\"ur LHS hinzuf\"ugen} They are defined as follows: Given a sequence $\mathseq{X_n}{n=1}^N$ of random variables with $X_n\sim\uniform([0,1)^d)$, the Latin hypercube sample in $s$ dimensions is defined as
%\[
%Y_n=\frac{\pi(n)-1+X_n}{N},
%\]
%where $\pi$ is a permutation on $\{1,\ldots,N\}$ chosen at random.  In \cite{gnew:discr_neg_dep_latin_hypcube}, it was proven that for any $a,b\in[0,1)^s$, defining $I=[0,b)\setminus[0,a)$, the sequence $\mathseq{\mathds{1}_{I}(Y_n)}{n=1}^N$ is $\gamma$-negatively dependent. 

\section{Proofs of results}\label{fi_cp:star_disc}
In this section, we will apply the maximal Bernstein inequality from \RefThm{fi_thm:max_bernstein_ineq} to derive bounds for the discrepancy of double-infinite matrices of elements in $[0,1)$ with entries stemming from uniformly distributed, independent random variables. In order to achieve this aim, we amend the structure of \cite[Theorem 1]{aistl:prob_star_discr_double_inf_mat} to our situation. First, we introduce some notation which allows for a coherent presentation of the topic. For $k\geq -1$, let
\[
b_k=\begin{cases}
3 &\text{if } k = -1,0\\
2^{-k+3}(1-2^{-k}) &\text{if } k> 0
\end{cases}.
\]
For a sequence $\mathseq{X_{n,s}}{n,s\in\N}$ of random variables we define

\[
X^{(N)}_s=(X_{N,1},
\ldots 
,X_{N,s})
\]
and\[
\calP_{N,s}=\begin{pmatrix}
X^{(1)}_s \\
X^{(2)}_s\\
\vdots \\
X^{(N)}_s\end{pmatrix}
=\begin{pmatrix}
X_{1,1}&X_{1,2}&\cdots &X_{1,s}\\
X_{2,1}&X_{2,2}&\cdots &X_{2,s}\\
\vdots&\vdots&\ddots&\vdots\\
X_{N,1}&X_{N,2}&\cdots &X_{N,s}
\end{pmatrix}.
\]
Another technical tool which is needed for the proof are $\delta$-covers and $\delta$-bracketing covers that we introduce next.
\begin{definition}Let $\delta>0$ and $\calA\subseteq\mathcal{P}([0,1]^s)$. A $\delta$-cover for $\calA$ is a set $\Gamma\subseteq \calA$ such that for every $A\in\calA$ there exist $U,V\in \Gamma\cup\{\emptyset\}$ satisfying $U\subseteq A\subseteq V$ and
\[
\leb(V\setminus U)\leq \delta.
\]
A $\delta$-bracketing cover for $\calA$ is a set $\Delta\subseteq \calA\times\calA$ such that for every $(A,B)\in\Delta $ with $A\subseteq B$ the inequality
\[
\leb(B\setminus A)\leq \delta
\]
holds and for every $A\in\calA$ there exist $(U,V)\in \Delta$ with $U\subseteq A\subseteq V$. Furthermore, we define $N(\delta,\calA)$ to be the minimal size of a $\delta$-cover for $\calA$ and $N_{[~]}(\delta,\calA)$ as the minimal size of a $\delta$-bracketing cover for $\calA$.
\end{definition}
These two concepts are closely related. In fact, it is easy to show that for any $\delta>0$ we have
\[
 N(\delta,\calA)\leq 2 N_{[~]}(\delta,\calA)\leq N(\delta,\calA)(1+N(\delta,\calA)).
\]
In this paper we will mostly use two families of sets, namely
\begin{align*}
&\calC^s=\mathset{[0,x)}{x\in[0,1]^s}
&\calR^s=\mathset{[x,y)}{x,y\in[0,1]^s}
\end{align*}
The following recent inequality from \cite{gnew:gen_faul_ineq}, which improved an earlier result from \cite{gnew:brack_numbers}, give us an estimate on the minimal size of $\delta$-bracketing covers and therefore by extension of  $\delta$-covers.
%\begin{theorem}\cite[Theorem 1.15]{gnew:brack_numbers} For $s\in\N$ the minimal size of a $\delta$-bracketing cover of $\mathcal{}R^s$ has the upper bound
%\begin{equation}\label{fi_eq:ineq_size_brack_cover_old}
%N_{[~]}(\delta,\mathcal{}R^s)\leq 2^{s-1}\frac{s^s}{s!}(\delta^{-1}+1)^s.
%\end{equation}
%\end{theorem}
\begin{theorem}\cite[Theorem 2.5.]{gnew:gen_faul_ineq} \label{fi_thm:gpw} Let $C_s:=\max\left\{1,1.1^{s-101}\right\}$. Then for $s\in\N$ the minimal size of a $\delta$-bracketing cover for $\calC^s$ has the upper bound
\begin{equation}\label{fi_eq:ineq_size_brack_cover_new}
N_{[~]}(\delta,\calC^s)\leq C_s\frac{s^s}{s!}(\delta^{-1}+1)^s.
\end{equation}
\end{theorem}
Additionally, there are also bounds for the sets in $\calR^s$.
\begin{theorem}\cite[Lemma 1.17 and 1.18]{gnew:brack_numbers} \label{fi_thm:brack_num_not_star} 
\begin{equation}\label{fi_eq:ineq_size_cover_not_star}
N(\delta,\calR^s)\leq N\left(\frac{\delta}{2},\calC^s\right)^2
\end{equation}
\begin{equation}\label{fi_eq:ineq_size_brack_cover_not_star}
N_{[~]}(\delta,\calR^s)\leq N_{[~]}\left(\frac{\delta}{2},\calC^s\right)^2
\end{equation}
\end{theorem}
Finally, we define 
\[
\calD^s:=\mathset{A\setminus B}{A,B\in\calR^s}.
\]
This puts us into the position to prove \RefThm{fi_thm:discrep_bound_not_star}. We decided to include the complete proof for the (extreme) discrepancy in our paper because this case is to our knowledge not treated in the existing literature. For the case of the star-discrepancy we will keep our presentation very short to avoid unnecessary repetitions.
\begin{proof}[\RefThm{fi_thm:discrep_bound_not_star}]
Let $C_s$ be as in \RefThm{fi_thm:gpw}. First, we define the sequences
\begin{align} \label{fi:eq:def:ak}
a_{k,s}=b_k\cdot \left(\frac{\ln\left(\frac{C_s^2}{\pi}2^{k+2}\right)}{s}+2\ln\left(2e(2^{k+2}+1)\right)\right)
\end{align}
for $k\geq 0$ and $a_{-1,s}=a_{0,s}$. With that we define
\[
A_s=2\left(\sum_{k=-1}^\infty\sqrt{a_{k,s}}\right)^2, B=2\left(\sum_{k=-1}^\infty\sqrt{b_{k}}\right)^2<178.
\]
Now we start with the actual proof. By \RefEq{fi_eq:ineq_size_brack_cover_new} and \RefThm{fi_thm:brack_num_not_star}, for every $\varepsilon>0$ we have 
\[
N_{[~]}(\varepsilon,\calR^s)\leq \left(C_s\frac{s^s}{s!}\left(\frac{2}{\varepsilon}+1\right)^s\right)^2,
\]
which means that for every $k\geq -1$, for $\varepsilon=2^{-(k+1)}$ we have
\begin{align}
N_{[~]}(\varepsilon,\calR^s)&\leq  \left(\frac{s^s}{s!}C_s(2^{k+2}+1)^s\right)^2\nonumber\\
&\leq\frac{\exp(2s)}{2\pi s}C_s^2\left(2^{k+2}+1\right)^{2s}\nonumber\\
&\leq\frac{1}{2}\exp\left({\ln\left(\frac{C_s^2}{\pi}\right)+2\alpha s\ln\left(e(2^{k+2}+1)\right)}\right)\label{fi_eq:size_brack_cover_new}
\end{align}
where we used Stirling's formula. Furthermore,
\[
N(\varepsilon,\calR^s)\leq 2N_{[~]}(\varepsilon,\calR^s)\leq   \frac{1}{2}\exp\left({\ln\left(2\frac{C_s^2}{\pi}\right)+2\alpha s\ln\left(e(2^{k+2}+1)\right)}\right).
\]
We now define
\[
\Omega_{M,s}=\left\{\max_{2^M\leq N<2^{M+1}}N\cdot D_N^{s}(\calP_{N,s})>\sqrt{\alpha A_s+\beta B\frac{\ln (M)}{s}}\sqrt{s 2^M}\right\}.
\]
Furthermore, we put
\[
y_{k,s}={\alpha a_{k,s}+\beta b_k\frac{\ln (M)}{s}}~\text{ and }~t_{k,s}=\sqrt{y_{k,s}}\sqrt{s 2^M}
\]
and choose $L\in\N$ such that
\begin{equation}\label{fi_eq:y_0_leq_two_power_L}
\frac{1}{2}\left(\frac{1}{2\sqrt{2}}\sqrt{y_{-1,s}}\sqrt{\frac{s}{2^M}}\right)<2^{-L}\leq\frac{1}{2\sqrt{2}}\sqrt{y_{-1,s}}\sqrt{\frac{s}{2^M}}.
\end{equation}
The proof will now proceed as follows: 
Our goal is to show that 
\[
\Prob{\bigcup_{M=1}^\infty\bigcup_{s=1}^\infty\Omega_{M,s}}<1.
\]
First, in step 1, we will use a dyadic chaining in order to create some $\calA_k\subseteq \calD^s$ and some events $E_k(I)$, $I\in\calA_k$ with 
\[
\Omega_{M,s}\subseteq\bigcup_{k=0}^L\bigcup_{I\in\calA_k}E_k(I).
\]
Then, in step 2, we will use the maximal Bernstein inequality in order to find an estimate for $\Prob{E_k(I)}$.\\
After that, in step 3, we will realize that the definition of the $a_{k,s}$ and $b_k$ are exactly tailored to obtain
\[
2^{k+1}(1+s)^\alpha M^\beta\vert \calA_k\vert\Prob{E_k(I)}\leq 1,
\]
which implies
\[
\Prob{\Omega_{M,s}}\leq\frac{1}{(1+s)^\alpha M^\beta},
\]
from which we will finally deduce our result in step 4.\\%[12pt]
\textbf{Step 1:} We will use dyadic chaining. For $1\leq k<L$ let $\Gamma_k$ be a $2^{-k}$-cover of $\calR^s$. Furthermore, let $\Delta_L$ be a $2^{-L}$-bracketing cover of $\calR^s$ and define
\[
\Gamma_L=\mathset{A\in\calR^s}{\exists B\in\calR^s~(A,B)\in\Delta_L}
\]
\[
\Gamma_{L+1}=\mathset{B\in \calR^s}{\exists A\in\calR^s~(A,B)\in\Delta_L}.
\]
For each $A\in\calR^s$ we choose a pair $(P_L(A),P_{L+1}(A))\in \Gamma_L\times\Gamma_{L+1}$ such that $\leb(P_{L+1}(A)\setminus P_{L}(A))\leq 2^{-L}$ and $P_{L}(A)\subseteq A\subseteq P_{L+1}(A)$. Inductively for each $1\leq k<L$ we choose some $P_{k}(A)\in\Gamma_k\cup\{\emptyset\}$ such that $P_{k}(A)\subseteq  P_{k+1}(A)$ and $\leb(P_{k+1}(A)\setminus P_{k}(A))\leq 2^{-k}$. Finally, we define $P_{0}(A)=\emptyset$.\\
For each $A\in\calR^s$ we write $I_k(A)=P_{k+1}(A)\setminus P_{k}(A)$, which means that we have $\leb(I_k(A))\leq 2^{-k}$ and
\begin{equation}\label{fi_eq:subseteq_dyad_chain}
\bigcup_{k=0}^{L-1}I_k(A)\subseteq A\subseteq\bigcup_{k=0}^{L}I_k(A).
\end{equation}
For each $1\leq k\leq L$ let $\calA_k=\mathset{I_k(A)}{A\in\calR^s}$. 
By \RefEq{fi_eq:size_brack_cover_new} we can assume that 
\begin{equation}\label{fi_eq:card_Ak}
\card{\calA_k}\leq \card{\Gamma_{k+1}}\leq\frac{1}{2}\exp\left({\ln\left(\frac{C_s^2}{\pi}\right)+2\alpha s\ln\left(e(2^{k+2}+1)\right)}\right)
\end{equation}
for all $1\leq k\leq L$. 
Next we define the events
\[
E_k(I)=\left\{\max_{2^M\leq N<2^{M+1}}\left\vert\sum_{n=1}^N\mathds{1}_I\left(X_s^{(n)}\right)-N\leb(I)\right\vert>t_{k,s}\right\},
\]
where $I\in\calA_k$. $E_k=\bigcup_{I\in\calA_k}E_k(I)$ and $E=\bigcup_{k=0}^L E_k$. Implicitly the set $E$ still depends on $M$ and $s$. We first show $\Omega_{M,s}\subseteq E$. We know that for every $\omega\in E^C=\bigcap_{k=0}^L\bigcap_{I\in\calA_k}E_k(I)^C$ we have
\[
\max_{2^M\leq N<2^{M+1}}\left\vert\sum_{n=1}^N\mathds{1}_I\left(X_s^{(n)}(\omega)\right)-N\leb(I)\right\vert\leq t_{k,s}~\text{ for all }0\leq k\leq L,I\in\calA_k.
\]
Given some $N\in\left[2^M,2^{M+1}\right)$, that means that by \RefEq{fi_eq:subseteq_dyad_chain} we obtain 
\begin{align}
\sum_{n=1}^N\mathds{1}_{A}(X_s^{(n)}(\omega))&\leq\sum_{k=0}^L\sum_{n=1}^N\mathds{1}_{P_{k+1}(A)\setminus P_k(A)}(X_s^{(n)}(\omega))\nonumber\\
&\leq\sum_{k=0}^L\left(N\leb(P_{k+1}(A)\setminus P_k(A))+t_{k,s}\right)\nonumber\\
&=N\leb(A)+N\leb(P_{L+1}(A)\setminus A)+\sum_{k=0}^L t_{k,s}.	\nonumber
\end{align}
for every $A\in\calR^s$. Since $P_L(A)\subseteq A\subseteq P_{L+1}(A)$ we get
\[
\leb(P_{L+1}(A)\setminus A)\leq \leb(P_{L+1}(A)\setminus P_L(A))\leq 2^{-L}\leq\frac{1}{2}\sqrt{y_{-1,s}}\sqrt{\frac{s}{2^M}}.
\]
Now, since $N<2^{M+1}$, we have $N\leb(P_{L+1}(A)\setminus A)\leq \sqrt{y_{-1,s}}\sqrt{\frac{s}{2^M}}$ and therefore
\begin{align}
\sum_{n=1}^N\mathds{1}_{A}(X_s^{(n)}(\omega))&\leq N\leb(A)+\sqrt{s 2^M}\sum_{k=-1}^L t_{k,s}\nonumber\\
&\leq N\leb(A)+\sum_{k=-1}^L\left(\sqrt{\alpha a_{k,s}}+\sqrt{\beta b_{k}}\right)\nonumber\\
&\leq N\leb(A)+\sqrt{\left(\alpha \sum_{k=-1}^L a_{k,s}\right)^2}+\sqrt{\beta\left( \sum_{k=-1}^L b_{k}\right)^2}\nonumber\\
&\leq N\leb(A)+\sqrt{\alpha A+\beta B\frac{\ln(M)}{s}}\sqrt{s 2^M}.\nonumber
\end{align}
A similar argument shows
\begin{align}
\sum_{n=1}^N\mathds{1}_{A}(X_s^{(n)}(\omega))&\geq\sum_{k=0}^{L-1}\sum_{n=1}^N\mathds{1}_{P_{k+1}(A)\setminus P_k(A)}(X_s^{(n)}(\omega))\nonumber\\
&\geq N\leb(A)-N\leb(P_{L+1}(A)\setminus A)-\sum_{k=0}^{L-1} t_{k,s}.	\nonumber\\
&\geq N\leb(A)-\sqrt{\alpha A_s+\beta B\frac{\ln(M)}{s}}\sqrt{s 2^M}.\nonumber
\end{align}
In combination, this means
\[
N\left\vert \frac{1}{N}\sum_{n=1}^N\mathds{1}_{I}(X_s^{(n)}(\omega))-\leb(\overline{I}) \right\vert\leq \sqrt{\alpha A_s+\beta B\frac{\ln(M)}{s}}\sqrt{s 2^M}.
\]
Since $A$ was arbitrary, we get
\[
N\cdot D_N^{s}(X_S^{(n)}(\omega))\leq\sqrt{\alpha A_s+\beta B\frac{\ln(M)}{s}}\sqrt{s 2^M},
\]
and due to the fact that this holds for every $\omega\in E^C$, the inclusion $\Omega_{M,s}\subseteq E$ holds.\\%[12pt]
\textbf{Step 2:} Next, we are going to estimate $\Prob{E_k(I)}$ for arbitrary $I\in\calA_k$. Define $Z_n=\mathds{1}_I(X^{(n)})-\leb(I)$. Then we have $\Exp{Z_n}=0$, $\vert Z_n\vert\leq \max\{\leb(I),1-\leb(I)\}\leq 1$, and $\sigma^2=\Exp{Z_n^2}=\leb(I)(1-\leb(I))$. By using the maximal Bernstein inequality \RefEq{fi_eq:max_bernstein_ineq} with $t=t_{k,s}$ we get
\[
\Prob{E_k(I)}\leq 2 \exp\left({-\frac{t_{k,s}^2/2^M}{4\leb(I)(1-\leb(I))+\frac{2}{3}t_{k,s}/2^M}}\right).
\]
A short calculation involving \RefEq{fi:eq:def:ak} implies that $y_{k,s}$ is decreasing in $k$. Therefore, applying \RefEq{fi_eq:y_0_leq_two_power_L} yields
\[
\frac{t_{k,s}}{2^M}=\frac{\sqrt{y_{k,s}} \sqrt{s 2^M}}{2^M}\leq\sqrt{y_{0,s}}\frac{\sqrt{s}}{ \sqrt{2^M}}<\sqrt{2} \, 2^{-L+2}
\]
and thereby
\[
\frac{2}{3}\frac{t_{k,s}}{2^M}<\frac{2}{3}\sqrt{2}2^{-L+2}\leq\begin{cases}
2^{-k+1}&\text{for }k=0,\ldots,L-1\\
2^{-L+2}&\text{for }k=L
\end{cases}.
\]
Moreover
\[
\frac{t_{k,s}^2}{2^M}=\frac{\left(\sqrt{y_{k,s}}\sqrt{s2^M}\right)^2}{2^M}=y_{k,s}\cdot s.
\]
The next aim is to minimize the value of $4\leb(I)(1-\leb(I))+\frac23 t_{k,s}/2^M$. For this we have to distinguish between several cases.\\
First, if $k=0$, we recall that $\leb(I)(1-\leb(I))\leq\frac{1}{4}$ in order to see
\[
\Prob{E_k(I)}\leq 2\exp\left({-\frac{y_{k,s}\cdot s}{3}}\right).
\]
Now, if $k=1$ we get (again with $\leb(I)(1-\leb(I))\leq\frac{1}{4}$) that
\[
\Prob{E_k(I)}\leq \exp\left({-\frac{y_{k,s}\cdot s}{2}}\right)=\exp\left({-\frac{y_{k,s}\cdot s}{2^{-k+3}(1-2^{-k})}}\right).
\]
For $k\in\{2,\ldots,L-1\}$ we have $\leb(I)\leq 2^{-k}$ since $I\in\Gamma_k$. Basic calculus tells us that the functions $\lambda(1-\lambda)$ and $4\lambda(1-\lambda)+2^{-k+1}(1-\lambda)$ take their maxima on the interval $[0,2^{-k}]$ at $\lambda=2^{-k}$, and that maximum is $2^{-k}(1-2^{-k})$ or $3\cdot 2^{-k+1}(1-2^{-k})$ respectively. Furthermore, $4\lambda(1-\lambda)+2^{-k+1}\lambda$ also takes its maximum at $2^{-k}$, resulting in a maximum of $2\cdot 2^{-k+1}(1-2^{-k})+2^{-k+1}2^{-k}\leq 3\cdot 2^{-k+1}(1-2^{-k})$. This implies
\[
\Prob{E_k(I)}\leq \exp\left({-\frac{y_{k,s}\cdot s}{2^{-k+3}(1-2^{-k})}}\right).
\]
Lastly, we come to $k=L$. In a similar way as in the last case it can be shown that 
\[
\max_{\lambda\in[0,2^{-L}]}\{4\lambda(1-\lambda)+ 2^{-L+2}\}\leq 2^{-L+3}(1-2^{-L}).
\]
Combining all cases, we see that
\begin{equation}\label{fi_eq:ineq_prob_EkI}
\Prob{E_k(I)}\leq \exp\left({-\frac{y_{k,s}\cdot s}{b_k}}\right)
\end{equation}
always holds.\\%[12pt]
\textbf{Step 3:} Using \RefEqTwo{fi_eq:card_Ak}{fi_eq:ineq_prob_EkI} we calculate
\begin{align}
\vert \calA_k\vert\Prob{E_k(I)}
&\leq \frac{1}{2}\exp\left({\ln\left(2\frac{C_s^2}{\pi}\right)+2\alpha s\ln\left(e(2^{k+2}+1)\right)}\right)\exp\left({-\frac{s\alpha a_{k,s}+\beta b_k{\ln (M)}}{b_k}}\right)\nonumber\\
&=\exp\left(\ln\left(2\frac{C_s^2}{\pi}\right)+2\alpha s\ln\left(e(2^{k+2}+1)\right)-\frac{s\alpha a_{k,s}+\beta b_k{\ln (M)}}{b_k}\right)\nonumber\\
&= \exp\left({\ln\left(2\frac{C_s^2}{\pi}\right)+2\alpha s\left(\ln\left(e(2^{k+2}+1)\right)-\frac{a_{k,s}}{b_k}\right)-\beta\ln(M)}\right)\nonumber\\
&\leq \exp\left({\alpha s\left(\frac{\ln\left(2\frac{C_s^2}{\pi}\right)}{s}+2\ln\left(e(2^{k+2}+1)\right)-\frac{a_{k,s}}{b_k}\right)-\beta\ln(M)}\right).\nonumber
\end{align}
Together with
\[
2^{k+1}(1+s)^\alpha M^\beta\leq \exp\left({\ln(2^{k+1})+\alpha s\ln(2)+\beta\ln(M)}\right),
\]
this results in
\begin{align*}
2^{k+1}(1+s)^\alpha M^\beta&\vert \calA_k\vert\Prob{E_k(I)}\\
&\leq \exp\left({\alpha s\left(\frac{\ln\left(\frac{C_s^2}{\pi}2^{k+2}\right)}{s}+2\ln\left(2e(2^{k+2}+1)\right)-\frac{a_{k,s}}{b_k}\right)}\right)\\
&=\exp(0)=1.
\end{align*}
Hence 
\begin{equation}\label{fi_eq:ineq_Prob_EkI_alpha_beta}
\vert \calA_k\vert\Prob{E_k(I)}\leq 2^{-(k+1)}\frac{1}{(1+s)^\alpha} \frac{1}{M^\beta}.
\end{equation}
\textbf{Step 4:} Finally, from \RefEq{fi_eq:ineq_Prob_EkI_alpha_beta} we obtain
\begin{equation}\label{fi_eq:ineq_Prob_E_alpha_beta}
\Prob{E}\leq\sum_{k=0}^L\sum_{I\in\calA_k}\Prob{E_k(I)}\leq\sum_{k=0}^L\vert \calA_k\vert\Prob{E_k(I)}\leq \sum_{k=0}^L 2^{-(k+1)}\frac{1}{(1+s)^\alpha} \frac{1}{M^\beta}<\frac{1}{(1+s)^\alpha} \frac{1}{M^\beta}.
\end{equation}
Please note that $\Prob{\Omega_{1,s}}=0$. Thus,
\begin{align}
\Prob{M\cdot D_M^{s}(\calP_{M,s})> \sqrt{\alpha A_s+\beta B\frac{\ln (\log_2(M))}{s}}\sqrt{s M}}
&\leq\Prob{\bigcup_{s=1}^\infty\bigcup_{M=2}^\infty\Omega_{M,s}}\nonumber\\
&\leq\sum_{s=1}^\infty\sum_{M=2}^\infty\Prob{\Omega_{M,s}}\nonumber\\
&\leq \sum_{s=1}^\infty\sum_{M=2}^\infty\frac{1}{(1+s)^\alpha} \frac{1}{M^\beta}\nonumber\\
&=(\zeta(\alpha)-1)(\zeta(\beta)-1).\nonumber
\end{align}
In other words, in order to have
\[
\Prob{D_M^{s}(\calP_{M,s})\leq \sqrt{\alpha A_s+\beta B\frac{\ln (\log_2(M))}{s}}\sqrt{\frac{s}{M}}}>0
\]
we need to ensure $(\zeta(\alpha)-1)(\zeta(\beta)-1)<1$ as claimed. This completes the proof.
\end{proof}
\begin{proof}[\RefThm{fi_thm:discrep_bound_star}]
The proof essentially stays the same, except of the definition of $a_{k,s}$, which needs to be changed to
\[
a_{k,s}=b_k\cdot \left(\frac{\ln\left( C_s 2^{k+2}\sqrt{\frac{2}{\pi}}\right)}{s}+\ln\left(2e(2^{k+1}+1)\right)\right).
\]
and the fact that we use \RefThm{fi_thm:gpw} instead of \RefThm{fi_thm:brack_num_not_star}. Consequently, the maximal value of $A_s$ is changed as well.
\end{proof}
\begin{remark}
When comparing the proof of \RefThm{fi_thm:discrep_bound_not_star} with the one from \cite{aistl:prob_star_discr_double_inf_mat}, one observes that for $s>1$, we always have $A_s<A$ (precise values can be seen in \RefSec{fi_cp:num_res}). %Since discrepancy is very well understood in dimension $1$, this gives an improvement in all cases of interest.
\end{remark}
\begin{remark} \label{fi_rem:star}
The proof can easily be adapted for other classes of sets, provided that there is an upper bound for their bracketing numbers similar to \RefEq{fi_eq:ineq_size_brack_cover_new}. 
\end{remark}

\section{Numerical results}\label{fi_cp:num_res}
Finally, we will now present explicit calculations for the values of the involved parameters in the results of \RefSec{fi_cp:star_disc}. As a benchmark, we compare our results with the constants from \cite{aistl:prob_star_discr_double_inf_mat}. In the paper of Aistleitner and Weimar, the dependency of $A_{s}$ on $s$ was neglected although it was implicitly present also in their proof. However, in order to account for a fair comparison, we need to compare our results not just to the constants presented in \RefThm{fi_thm:aw}, but also to those which could have been obtained by not neglecting said dependency. Obviously there are different choices for $\alpha$ and $\beta$ that ensure 
$$(\zeta(\alpha)-1)(\zeta(\beta)-1)<1$$
and thus by \RefThm{fi_thm:discrep_bound_star} and \RefThm{fi_thm:discrep_bound_not_star} the existence of double-infinite sequences with the desired star-discrepancy. Here we will restrict ourselves to ${\alpha=\beta=1.73>\zeta^{-1}(2)}$ because the symmetry simplifies calculations and was also considered in \RefThm{fi_thm:aw}. Not ignoring the dependency of  $A_{s}$ in \cite{aistl:prob_star_discr_double_inf_mat} would result in 
\begin{equation}\label{fi_eq:aks_old_dim}
a_{k,s}=b_k\cdot \left(\frac{\ln\left(2^{k+1}\sqrt{\frac{2}{\pi}}\right)}{s}+\ln\left(4e(2^{k+1}+1)\right)\right)
\end{equation}
and $A_s=2\left(\sum_{k=-1}^\infty\sqrt{a_{k,s}}\right)^2$.

\begin{figure}[!ht]
     \centering
         \centering\includegraphics[width=0.8\textwidth]{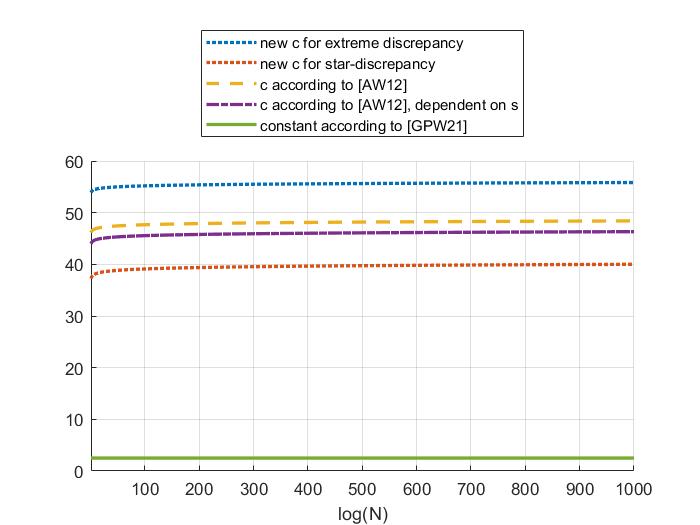}
         \caption{Values for $c$ dependent on $N$ for $s=10$.}
         \label{fi_fig::const_c_dep_N}
\end{figure}
In \RefFig{fi_fig::const_c_dep_N} we display numerical estimates for the smallest possible value $c$ such that 
\[
\Prob{D_N^{*s}(\calP_{N,s})\leq c\sqrt{\frac{s}{N}}}>0
\]
dependent on the sample size $N$. Besides the values for the discrepancy and star-discrepancy implied by our result, here we include the values from \cite{aistl:prob_star_discr_double_inf_mat} as a point of comparison, both with and without dependency from $s$. Finally, we include in the figure the result from \cite{gnew:gen_faul_ineq}, where $c=2.4968$ was proved for finite-dimensional point sets. It can be observed that while our result (in the one-dimensional case) is lower than that from \cite{aistl:prob_star_discr_double_inf_mat}, it is still much higher than the one from \cite{gnew:gen_faul_ineq}.
%We show our values for both the discrepancy and the star-discrepancy. Again we use the the values from \cite{aistl:prob_star_discr_double_inf_mat}, both with and without dependency from $s$, as a point of comparison.
\begin{figure}[!ht]
     \centering
         \includegraphics[width=0.8\textwidth]{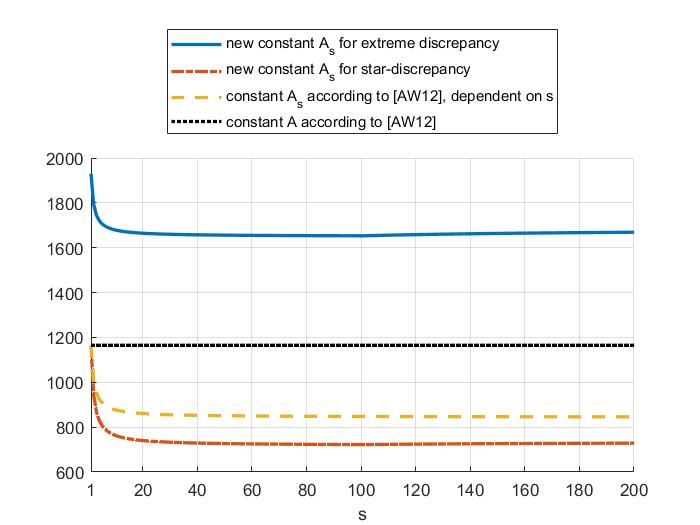}
         \caption{Overview of $A_s$ dependent on $s$}
         \label{fi_fig::const_A_dep_s}
\end{figure}
\begin{figure}[!ht]
     \centering
         \includegraphics[width=0.8\textwidth]{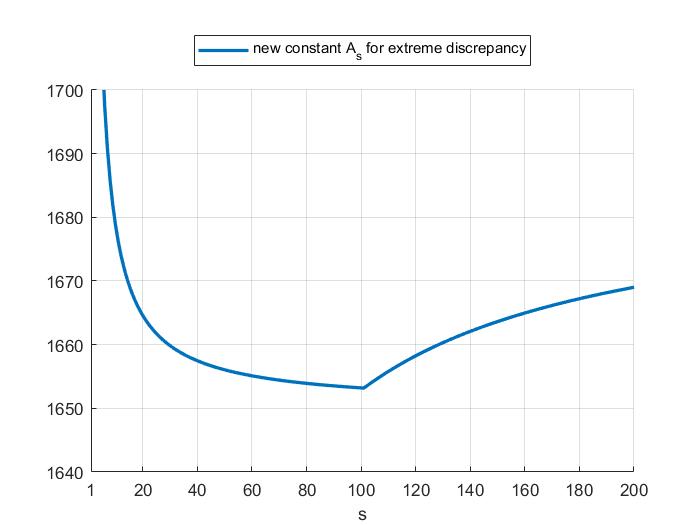}
         \caption{Close up view of $A_s$ dependent on $s$}
         \label{fi_fig::const_A_dep_s_big}
\end{figure}
In \RefFig{fi_fig::const_A_dep_s}, we show the four different values for $A_s$. For $s=1$ the values of the three parameters for the star-discrepancy are equal. Since we used
\[
a_{k,s}=b_k\cdot \left(\frac{\ln\left( C_s 2^{k+2}\sqrt{\frac{2}{\pi}}\right)}{s}+\ln\left(2e(2^{k+1}+1)\right)\right)
\]
instead of \RefEq{fi_eq:aks_old_dim}, this however changes for $s>1$. Furthermore, we see that the value of $A_s$ for the discrepancy is larger than that for the star-discrepancy as expected by \RefThm{fi_thm:brack_num_not_star}.\\%[12pt]
One aspect that can be observed in both figures is that the values for the parameters for the discrepancy are approximately twice those for the star-discrepancy. This stems from the fact that according to \RefThm{fi_thm:brack_num_not_star} it holds that $N(\delta,\calR^s)\leq N\left(\frac{\delta}{2},\calR^s\right)^2$ and by taking the logarithm in the proof of \RefThm{fi_thm:discrep_bound_not_star} and \RefThm{fi_thm:discrep_bound_star}, this basically doubles the value of $A_s$. \\%[12pt]
Finally in \RefFig{fi_fig::const_A_dep_s_big}, a close-up view of the discrepancy case is presented. It can be observed that up until $s=101$ the value for $A_s$ is decreasing and afterwards it is increasing until it reaches a value of approximately $1670$. This is due to the nature of $C_s$, compare \RefThm{fi_thm:gpw}, which increases from $s=101$ onward.\\%[12pt]
Summing up, all figures show that at least for higher dimensions, the numerical results we get are a clear improvement in comparison to \cite{aistl:prob_star_discr_double_inf_mat}. Still, in the one-dimensional situation as considered in \cite{gnew:gen_faul_ineq}, much more is known and therefore the values of the constants are significantly smaller. Nevertheless, we believe that our work is a step towards improving numerical values for higher dimensions.%, but also in broadening the methodology by including non-independent random variables.

\section{Future Research}
We are currently working on generalizing the approach used in this paper from i.i.d. sequences towards other types of random variables, for example negatively associated and  $\gamma$-negatively dependent random variables. 
The primary challenge when trying to do so, lies in finding appropriate maximal inequalities for these types of random variables, for example versions of the maximal Bernstein inequality.
There has been extensive research on this topic, see e.g. \cite{bouk:max_ineq_exp_decay,kevei:note_max_bern_ineq, 
moricz:exp_est_max_part_sums,szewczak:max_levi_otta_ineq_sums_vector}. We believe that methods from the existing literature can be adapted to the mentioned classes of examples. The next challenge in this context is then to prove the existence of infinite sequences or matrices with such characteristics. To the best of our knowledge, there are no known examples of this apart from some relatively simple ones, so this is another question which is currently wide open.

%We are currently working towards expanding the approach from this work and the one that came before it from iid sequences towards other types of random variables, for example negatively associated and  $\gamma$-negatively dependent random variables.\\
\newpage
\bibliographystyle{alpha}
\bibdata{references}
\bibliography{references}
\end{document}